\documentclass[11 pt, twoside]{amsart}
\usepackage{parskip}
\usepackage{amsthm,amsmath,amssymb}
\usepackage[a4paper, margin=2.5cm]{geometry}
\usepackage{mathrsfs}
\usepackage{verbatim}
\usepackage[toc,page]{appendix}
\usepackage{thumbpdf}
\usepackage{hyperref}
\newcommand\nnfootnote[1]{%
  \begin{NoHyper}
  \renewcommand\thefootnote{}\footnote{#1}%
  \addtocounter{footnote}{-1}%
  \end{NoHyper}
}
\theoremstyle{plain}
\newtheorem{thm}{Theorem}[section]

\newtheorem{mainthm}[]{Main Theorem}[]
\newtheorem{lem}[thm]{Lemma}
\newtheorem{prop}[thm]{Proposition}

\theoremstyle{definition}
 
\theoremstyle{definition}
\newtheorem{rem}[]{Remark}[]
\newtheorem*{rem*}{Remark}

\let\Im\relax

\let\H\relax
\let\Z\relax

\let\dh\relax

\def\Im{{\rm Im}}

\def\deg{{\rm deg}}

\def\Gr{{\rm Gr}}
\def\Sym{{\rm Sym}}
\def\hyp{{\rm hyp}}

\def\dim{{\rm dim}}

\def\FS{{\rm FS}}
\def\log{{\rm log}}
\def\inf{{\rm inf}}
\def\cosh{{\rm cosh}}
\def\sinh{{\rm sinh}}
\newcommand{\ncom}{\newcommand}
\ncom{\ep}{\epsilon}
\ncom{\rar}{\rightarrow}
\ncom{\thrar}{\twoheadrightarrow}
\ncom{\lrar}{\longrightarrow}
\ncom{\ov}{\overline}
\ncom{\what}{\widehat}
\let\H\relax

\let\dh\relax

\DeclareMathOperator{\xbar}{{\it{\overline{X}}}}
\DeclareMathOperator{\dh}{{\it{d}}_{\mathrm{hyp}}}

\DeclareMathOperator{\bk} {\mathcal{B}_{{\Gamma}}^{2{\it{k}}}}
\DeclareMathOperator{\bkxbar} {\mathcal{B}_{\Omega_{\it{\overline{X}}}}^{\it{k}}}

\DeclareMathOperator{\bkl} {\mathcal{B}_{\mathcal{L}}^{\it{k}}}
\DeclareMathOperator{\bkld} {\mathcal{B}_{\mathcal{L},{\it{D}}}^{\it{k}}}
\DeclareMathOperator{\hypx}{\mu_{{\it X}}^{\mathrm{hyp}}} 
\DeclareMathOperator{\hypbarx}{\mu_{{\it \overline{X}}}^{\mathrm{hyp}}} 
 
\DeclareMathOperator{\muberkbarx}{\mu_{{\it \overline{X}}}^{\mathrm{Ber},{\it{k}}}}

\DeclareMathOperator{\hypxdbar}{\mu_{{\it \overline{X}^d}}^{\mathrm{hyp}}} 
\DeclareMathOperator{\hypxdvolbar}{\mu_{{\it \overline{X}^d}, \mathrm{vol}}^{\mathrm{hyp}}}
\DeclareMathOperator{\H}{\mathbb{H}}

\DeclareMathOperator{\G}{\Gamma}

\DeclareMathOperator{\C}{\mathbb{C}}

\DeclareMathOperator{\R}{\mathbb{R}} 
 
\DeclareMathOperator{\symxd}{\mathrm{Sym}^{{\it d}}({\it X})} 
\DeclareMathOperator{\symxdbar}{\mathrm{Sym}^{{\it d}}({\it \overline{X}})}

\DeclareMathOperator{\Ocpt}{\Omega_{{\it{\overline{X}}}}} 
\DeclareMathOperator{\Ocptk}{\Omega_{{\it{\overline{X}}}}^{\otimes{\it{k}}}} 
 
\DeclareMathOperator{\Sk}{\mathcal{S}^{2{\it{k}}}(\Gamma)}
\DeclareMathOperator{\rx}{{\it{r}}_{\G}} 
\newcommand{\ignore}[1]{}
\ncom{\m}{\mbox}
\ncom{\sta}{\stackrel}
\ncom{\A}{{\mathbb A}}
\ncom{\Z}{{\mathbb Z}}
\ncom{\Q}{{\mathbb Q}}
\ncom{\HH}{{\mathbb H}}
\ncom{\al}{\alpha}
\ncom{\p}{{\mathbb P}}
\ncom{\K}{{\mathbb K}}
\ncom{\X}{{\mathbb X}}
\ncom{\f}{\frac}
\ncom{\cA}{{\mathcal A}}
\ncom{\cB}{{\mathcal B}}
\ncom{\cD}{{\mathcal D}}
\ncom{\cDB}{{\mathcal D \mathcal B}}
\ncom{\cX}{{\mathcal X}}
\ncom{\cO}{{\mathcal O}}
\ncom{\cW}{{\mathcal W}}
\ncom{\cL}{{\mathcal L}}
\ncom{\cP}{{\mathcal P}}
\ncom{\cH}{{\mathcal H}}
\ncom{\cS}{{\mathcal S}}
\ncom{\cM}{{\mathcal M}}
\ncom{\cC}{{\mathcal C}}
\ncom{\cT}{{\mathcal T}}
\ncom{\cF}{{\mathcal F}}
\ncom{\cN}{{\mathcal N}}
\ncom{\cJ}{{\mathcal J}}
\ncom{\cV}{{\mathcal V}}
\ncom{\cZ}{{\mathcal Z}}
\ncom{\cU}{{\mathcal U}}
\ncom{\cSU}{{\mathcal S \mathcal U}}
\ncom{\cG}{{\mathcal G}}
\ncom{\cQ}{{\mathcal Q}}
\ncom{\cR}{{\mathcal R}}
\ncom{\cY}{{\mathcal Y}}
\ncom{\cE}{{\mathcal E}}
\ncom{\cI}{{\mathcal I}}
\ncom{\mylabel}[1]{{\rm (#1)}\label{#1}}
\ncom{\Hom}{{\textit{Hom}}}
\ncom{\eop}{{\hfill $\Box$}}
\begin{document}
\baselineskip=16pt
\nnfootnote{Mathematics Subject Classification:  32A25, 11F11, 32N05, 53C07.}
\nnfootnote{Keywords:  Bergman kernel, hyperbolic Riemann surface, Fubini-Study metric, symmetric product.}
\setcounter{tocdepth}{1}
\title[Estimates of Bergman metric]{Estimates of K\"ahler metrics on noncompact finite volume hyperbolic Riemann surfaces, and their symmetric products}
\author{Anilatmaja Aryasomayajula}
\address{Department of Mathematics, Indian Institute of Science Education and Research (IISER) Tirupati, 
Transit campus at Sri Rama Engineering College, Karkambadi Road,
Mangalam (B.O),Tirupati-517507, India.}
\email{anilatmaja@gmail.com}
\author{Arijit Mukherjee}
\address{Department of Mathematics, Indian Institute of Science Education and Research (IISER) Tirupati, 
Transit campus at Sri Rama Engineering College, Karkambadi Road,
Mangalam (B.O),Tirupati-517507, India.}
\email{arijitmukherjee@iisertirupati.ac.in}
\begin{abstract}
Let $X$ denote a noncompact finite volume hyperbolic Riemann surface of genus $g\geq 2$, with only one puncture at $i\infty$ (identifying $X$ with its universal cover $\H$). Let $\xbar:=X\cup\lbrace i\infty\rbrace$ denote the Satake compactification of $X$. Let $\Ocpt$ denote the cotangent bundle on $\xbar$. For $k\gg1$, we derive an estimate for $\muberkbarx$, the Bergman metric associated to the line bundle $\cL^{k}:=\Ocptk\otimes \cO_{\xbar}((k-1)i\infty)$. 
 
 \vspace{0.1cm}\noindent
For a given $d\geq 1$, the pull-back of the Fubini-Study metric on the Grassmannian, which we denote by $\mu_{\symxdbar}^{\FS,k}$, defines a K\"ahler metric on $\symxdbar$, the $d$-fold symmetric product of $\xbar$. Using our estimates of $\muberkbarx$, as an application, we derive an estimate for $\mu_{\symxdbar,\mathrm{vol}}^{\FS,k}$, the volume form associated to the (1,1)-form $\mu_{\symxdbar}^{\FS,k}$.
\end{abstract}
\maketitle
\section{Introduction}\label{section1}
\subsection{History of the problem, and statements of the main results}\label{subsec1.1}
Comparison of various K\"ahler structures defined on a complex manifold is an area of fundamental importance in the field of complex geometry. Furthermore, estimates of Bergman kernels associated to high tensor-powers of holomorphic line bundles defined over a complex manifold, and estimates of the associated Bergman metric, is a closely related problem, and is also of immense interest in the field of complex geometry. Tian, Zelditch, Demailly, Catlin,  Ma, and Marinescu have contributed significantly to the field. 

\vspace{0.1cm}
Estimates of Bergman kernels associated to high tensor-powers of holomorphic line bundles, defined over a compact complex manifold are well understood, and have been extensively studied. However, results in the noncompact setting are scarce. In the recent past, optimal estimates of Bergman kernels associated to high tensor-powers of holomorphic line bundles, defined over complete sympletic manifolds have been derived in \cite{ma1}, \cite{kord},  \cite{hsiao}, and \cite{hsiao1}. Furthermore, we refer the reader to \cite{auvray} and \cite{auvray1}, where estimates of Bergman kernels associated to high tensor-powers of holomorphic line bundles, defined over noncompact hyperbolic Riemann surfaces,  were derived, and which are more relevant to the setting of the current article. Lastly, we refer the reader to \cite{ma}, which is a treatise on the topic of estimates of Bergman kernels. 

\vspace{0.1cm}
Estimates of Bergman metric associated to Bergman kernels associated to high tensor-powers of holomorphic line bundles, defined over a compact complex manifold, are also difficult to derive, and even more in the setting of noncompact complex manifolds. 

\vspace{0.1cm}
In \cite{ab}, the authors computed estimates of the Bergman metric associated to high tensor-powers of the cotangent bundle defined over a compact hyperbolic Riemann surface. As an application, the authors derived estimates of a particular K\"ahler metric on the symmetric product of a compact hyperbolic Riemann surface.  

\vspace{0.1cm}
Symmetric products of Riemann surfaces can be realized as the moduli space of vortices, and are of immense interest both in algebraic geometry, and theoretical physics, and have been vastly investigated (see \cite{arbarello}, \cite{kempf}, \cite{manton}, and \cite{biswas}). 

\vspace{0.1cm}
The symmetric product of a compact hyperbolic Riemann surface can be embedded into a Grassmannian, via the cotangent bundle. The pull-back of the Fubini-Study metric defines a K\"ahler metric on the symmetric product of the compact hyperbolic Riemann surface. Using estimates of the Bergman metric associated to Bergman kernels associated to high tensor-powers of holomorphic line bundles, the authors of \cite{ab}, derived estimates of the pull-back of the Fubini-Study metric. Similar results on K\"ahler metrics were realized in \cite{abms}.

\vspace{0.1cm}
In this article, we extend the estimates from \cite{ab} to noncompact finite volume hyperbolic Riemann surfaces, which we now describe.

\vspace{0.1cm}
Let $X$ denote a noncompact finite volume hyperbolic Riemann surface of genus $g\geq 2$. Let $\hypx$ denote the natural metric on $X$, which is compatible with its complex structure, and with constant curvature $-1$. By Riemann uniformization theorem, $X$ can be realized as $\G\backslash \H$, where $\G\subset \mathrm{PSL}_{2}(\R)$ is a cofinite Fuchsian subgroup. Without loss of generality, we assume that $i\infty$ is the only puncture of $X$ (identifying $X$ with its universal cover $\H$). Furthermore, we assume that $\G_{i\infty}$, the stabilizer of the cusp $i\infty$, is of the form 
\begin{align*}
\G_{i\infty}=\bigg\lbrace\bigg(\begin{matrix} 1&n\\0&1 \end{matrix}\bigg)\bigg|\,n\in\Z\bigg\rbrace.
\end{align*}

\vspace{0.1cm}
Let $\xbar:=X\cup\lbrace i\infty\rbrace$, denote the Satake compactification of $X$. The hyperbolic metric $\hypx$ extends to a singular metric on $\xbar$, which we denote by $\hypbarx$. The metric $\hypbarx$ is smooth at all $z\in \xbar\backslash \lbrace i\infty\rbrace$, and admits an integrable singularity at the puncture $i\infty$.

\vspace{0.1cm}
Let $\Ocpt$ denote the cotangent bundle on $\xbar$, and let $\cL^{k}:= \Ocptk\otimes \cO_{\xbar}\big((k-1)i\infty\big)$. For $k\geq 1$, let $\bkl$ denote the Bergman kernel associated to the line bundle $\cL^{k}$. Furthermore, let $\|\cdot\|_{\hyp}$ denote the point-wise metric on $H^{0}\big(\xbar,\cL^{k}\big)$, which is induced by the hyperbolic metric. Then, for $z\in \xbar\backslash\lbrace i\infty\rbrace$, the associated Bergman metric is defined as
\begin{align*}
\muberkbarx(z):=-\frac{i}{2\pi}\partial_{z}\partial_{\overline{z}}\log\|\bkl(z)\|_{\hyp}. 
\end{align*}

\vspace{0.1cm}
Let $\symxdbar$ denote the symmetric product of $\xbar$. The hyperbolic metric $\hypbarx$ defines a singular metric $\hypxdbar$ on $\symxdbar$, and let 
$\hypxdvolbar$ denote the associated volume form.  

\vspace{0.1cm}
Let $D$ be an effective divisor on $\xbar$ of degree $d\geq 1$. For $k\geq 1$ and large enough, we have the short exact sequence of sheaves over $\overline{X}$
\begin{align*}
0 \longrightarrow  \cL^{k}\otimes\cO_{\xbar}(-D) \longrightarrow \cL^{k}\longrightarrow \cL^{k}\vert_D  \longrightarrow 0,
\end{align*}
which induces the following injective homomorphism
\begin{align}\label{hom}
H^0\big(\xbar, \cL^{k}\otimes\cO_{\xbar}(-D)\big)\longrightarrow H^0\big(\xbar,\cL^{k}\big).
\end{align}

For a given $d\geq 1$, and $k\geq 1$ be large enough such that $n_{k}:=(2k-1)(g-1)+k-1>d$, let $\mathrm{Gr}(r_{k},n_{k})$ denote the Grassmannian parametrizing $r_{k}$-dimensional vector subspaces of $\C^{n_{k}}$, where $r_k=n_k-d$. Then, from the homomorphism as in equation \eqref{hom}, we have the holomorphic embedding 
\begin{align*}
\varphi_{\cL}^{k,d}: \symxdbar\hookrightarrow &\mathrm{Gr}(r_{k},n_{k})\\[0.1cm] (z_1,\ldots,z_d)\mapsto &H^0\big(\xbar,\cL^{k}\otimes {\mathcal O}_{\xbar}(-D)\big)\subset H^0\big(\xbar, \cL^{k}\big),
\end{align*}
where $D$ denotes the divisor $x_1+\cdots +x_d$ on $\xbar$, and $z_j$ is the complex coordinate of the point $x_j\in\xbar$, for each $1\leq j\leq d$. 

\vspace{0.1cm}
For further details regarding the above holomorphic embedding, we refer the reader to \cite{ab}. For convenience of the reader, we recall relevant details in section \ref{subsec-3.1} of this article. 

\vspace{0.1cm}
The Fubini-Study metric is the natural metric on $\mathrm{Gr}(r_{k},n_{k})$, and let $\mu_{\symxdbar,\mathrm{vol}}^{\mathrm{FS},k}(z)$ denote the volume form, associated to the pull-back of the Fubini-Study metric on $\mathrm{Gr}(r_{k},n_{k})$, via the holomorphic embedding $\varphi_{\cL}^{k,d}$. 

\vspace{0.1cm}
\subsection{Statements of main results}\label{subsec1.2}
Since the Bergman kernel $\bkxbar$ vanishes at $i\infty$, the Bergman metric $\muberkbarx$ is not defined at $i\infty$. However, as shown in Theorem \ref{thm10}, the ratio 
\begin{align*}
\bigg|\frac{\muberkbarx(z)}{\hypbarx(z)}\bigg|
\end{align*}
is well defined.

\vspace{0.1cm}
Before we state the main theorem, we explain the following notation. For two functions $f,g$ defined on $\xbar$, the symbol $f\ll_{\G}g$ means 
$f\leq c g$, where the constant $c$ depends only on $\G$ (equivalently on $\xbar$). Similarly, $f\gg_{\G}g$ means $f\geq cg$, where the constant $c$ depends only on $\G$ (equivalently on $\xbar$).

\vspace{0.1cm}
We now state the first main result of the article, which is proved as Theorem \ref{thm10} in section \ref{subsec-2.3}.

\vspace{0.1cm}
\begin{mainthm}\label{mainthm1}
With notation as above, for $k\gg1$, we have the following estimate
\begin{align}\label{est:mainthm1}
\sup_{z\in \xbar}\bigg|\frac{\muberkbarx(z)}{\hypbarx(z)}\bigg|=O_{\G}(k^{3}),
\end{align}
where the implied constant depends only on the Fuchsian subgroup $\G$.
\end{mainthm}

\vspace{0.1cm}
We now state the second main result of the article, which is proved as Theorem \ref{thm11} in section \ref{subsec-3.2}.

\vspace{0.1cm}
\begin{mainthm}\label{mainthm2}
With notation as above, for a given $d\geq 1$, and $k\gg1$, we have the following estimate
\begin{align}\label{est:mainthm2}
\sup_{z\in \symxdbar}\bigg|\frac{\mu_{\Sym(\xbar),\mathrm{vol}}^{\mathrm{FS},k}(z)}{\mu_{\xbar^d,\mathrm{vol}}(z)}\bigg|=O_{\G}(k^{3d}),
\end{align}
where the implied constant depends only on the Fuchsian subgroup $\G$.
\end{mainthm} 

\begin{rem}
The main technique involved in the proofs of the above theorems, is the analysis of the Poincar\'e series representing the Bergman kernel $\bkl$. We utilize the fact that the complex vector space $H^{0}\big(\xbar,\cL^{k}\big)$, the space of global sections of the line bundle $\cL^{k}$, is isometric to $\Sk$, the complex vector space of cusp forms of weight-$2k$, with respect to the Fuchsian subgroup $\G$.  

Lastly, our estimates \eqref{est:mainthm1} and \eqref{est:mainthm2} gain a factor of $k$ and $k^{d}$, respectively, when compared to the estimates derived in \cite{ab}, which is in accordance with the gain of a factor of $\sqrt{k}$ for the estimate of the Bergman kernel for noncompact hyperbolic Riemann surfaces, when compared to the compact setting. 
\end{rem}
\section{Estimates of the Bergman metric on noncompact hyperbolic Riemann surface}\label{sec2}
\subsection{Notation and background material}\label{subsec-2.1}
In this section, we setup the notation, and recall the background material required to prove Main Theorem \ref{mainthm1}.

\vspace{0.1cm}
Let $X$ be a noncompact finite volume hyperbolic Riemann surface of genus $g\geq 2$. From uniformization theorem from complex analysis, $X$ is isometric to $\G\backslash\H$, where $\G\subset\mathrm{PSL}_{2}(\R)$ is a cofinite Fuchsian subgroup, which acts on $\H$ via fractional linear transformations. Locally, we identify $X$ with its universal cover $\mathbb{H}$.  Let $\cF_{\G}$ denote a fixed fundamental domain of $X$. Without loss of generality, we assume that
\begin{align*}
\cF_{i\infty}:=\big\lbrace z=x+iy\in\H\big|\,0\leq x\leq 1,\,y\geq 1\slash 2  \big\rbrace\subset \cF_{\G}.
\end{align*}

The punctures of $X$ are in one-one correspondence with cusps of $\G$. Without loss of generality, we assume that $i\infty$ is the only puncture of $X$, and we further assume that 
\begin{align*}
\G_{i\infty}=\bigg\lbrace\bigg(\begin{matrix} 1&n\\0&1 \end{matrix}\bigg)\bigg|\,n\in\Z\bigg\rbrace.
\end{align*}

Let $\xbar:=X\cup\lbrace i\infty\rbrace$ denote the Satake compactification of $X$. We now describe the local coordinates of $\xbar$. For $z\in \xbar\backslash\lbrace i\infty\rbrace$, let $U_{r}(z)$ denote a coordinate disc around $z$, and of radius $r>0$. For $w\in U_{r}(z)$, the local coordinate function $\vartheta_{z}(w)$ is given by the following formula
\begin{align*}
\vartheta_{z}(w):=w-z. 
\end{align*}  

Let $U_{r}(i\infty)$ denote a coordinate disc around $i\infty$, and of radius $r>0$. For $w\in U_{r}(i\infty)$, the local coordinate function $\vartheta_{i\infty}(w)$, also denoted by $q(z)$, is given by the following formula
\begin{align}\label{qcor}
q(z)=\vartheta_{i\infty}(w):=e^{2\pi iw}. 
\end{align}  

\vspace{0.1cm}
Let $\hypx$ denote the hyperbolic metric on $\H$, and for $z:=x+iy\in \H$, the hyperbolic metric is given by the following formula
\begin{align*}
\hypx(z):=\frac{i}{2}\cdot\frac{dz d\overline{z}}{y^2}=\frac{dxdy}{y^{2}}.
\end{align*}

Let $\dh(z,w)$ denote the hyperbolic distance between the points $z$ and $w$ on $\H$. For $z=x+iy,\,w=u+iv\in \H$, the hyperbolic distance is given by the following formula
\begin{align}\label{dh}
\cosh^{2}\big(\dh(z,w)\slash 2\big)=\frac{\big|z-\overline{w}\big|^{2}}{4yv}.
\end{align} 

The hyperbolic metric descends to define a K\"ahler metric on $X$. Locally, for any two points $z,w\in X$, the geodesic distance between the points $z$ and $w$ is given by $\dh(z,w)$.

\vspace{0.1cm}
We now define the injectivity radius $\rx$ of $X$, which is by the following formula
\begin{align*}
\rx:=\inf \big\lbrace\dh(z,\gamma z)\big|\, z\in\cF_{\G},\, \gamma\in \Gamma \backslash \Gamma_{i\infty}\big\rbrace.
\end{align*}  

The hyperbolic metric $\hypx$ extends to a singular metric on $\xbar$, which we denote by $\hypbarx$. Locally, for $z\in U_{r}(i\infty)$, the hyperbolic metric is given by the following formula
\begin{align*}
\hypbarx(z)=\frac{i}{2}\cdot\frac{dq(z)d\overline{q(z)}}{\big(\big|q(z)\big|\log\big|q(z)\big|\big)^{2}},
\end{align*}
where  $U_{r}(i\infty)$ denote a coordinate disc around $i\infty$, and of radius $r>0$, and $q(z)$ is as in equation \eqref{qcor}.

\vspace{0.2cm}
\paragraph{\bf{Bergman kernel and Cusp forms}}
Let $\Ocpt$ denote the cotangent bundle of holomorphic differential $1$-forms on $\xbar$. For $k\geq 1$, set
\begin{align*}
\cL^{k}:=\Ocptk\otimes \cO_{\xbar}\big((k-1)i\infty\big).
\end{align*}

Let $H^0\big(\xbar,\cL^{k}\big)$ denote the space of global sections of $\cL^{k}$, meromorphic differentials with a pole of order at most $(k-1)$ at the puncture $i\infty$, and holomorphic at all $z\in\xbar\backslash\lbrace i\infty\rbrace$.  Let $q(z)$ be as in equation \eqref{qcor}.  Then, for $k\geq 1$, at $z\in\xbar\backslash\lbrace i\infty\rbrace$, $\omega\in H^0\big(\xbar,\cL^{k}\big)$ is of the form 
\begin{align*}
\omega(z)=f(z)dz^{\otimes k}=\frac{f(q(z))dq(z)^{\otimes k}}{(2\pi i q(z))^{k}},
\end{align*} 

The differential $dq(z)^{\otimes k}\slash q(z)^k$ has a pole of order $k$ at $i\infty$. Since, $\omega$ has a pole of order at most $k-1$ at $i\infty$, it implies that $f\in \Sk$, the complex vector space of cusp forms of weight-$k$, with respect to $\G$. Hence, we have an isometry 
\begin{align}\label{isometry}
H^0\big(\xbar,\cL^{k}\big)\simeq \Sk.
\end{align}

\vspace{0.1cm}
For $w\in H^0\big(\xbar,\cL^{k}\big)$, where $\omega(z)=f(z)dz^{\otimes k}$, at $z=x+iy\in \xbar$ (identifying $X$ with its universal cover $\H$), the point-wise hyperbolic metric $\|\cdot\|_{\hyp}$ on $H^0\big(X,\cL^k\big)$, is given by the following formula 
\begin{align*}
\| \omega(z) \|_{\hyp}^{2}=y^{2k}\big|f(z)\big|^{2}.
\end{align*}

The point-wise metric $\|\cdot\|_{\hyp}$ induces an $L^{2}$-metric on $ H^0\big(\xbar,\cL^{k}\big)$. For $z=x+iy\in X$ (identifying $X$ with its universal cover $\H$), $\omega, \eta\in H^0\big(\xbar,\cL^{k}\big)$, where $\omega(z)=f(z)dz^{\otimes k}$, $\eta(z)=g(z)dz^{\otimes k}$, the hyperbolic $L^2$-metric $\langle\cdot,\cdot\rangle_{\hyp}$ on $H^0\big(\xbar,\cL^k\big)$, is given by the following formula
\begin{align}\label{L^2-metric}
\big\langle \omega,\eta\big\rangle_{\hyp}:=\int_{\cF_{\G}}y^{2k}f(z)\overline{g(z)}\hypbarx(z).
\end{align} 

\vspace{0.1cm}
Let $n_k$ denote the complex dimension of $H^0\big(\xbar,\cL^{k}\big)$. Let $\big\lbrace\omega_{1},\ldots, \omega_{n_k}\big\rbrace$ be an orthonormal basis of $H^0\big(\xbar,\cL^{k}\big)$ with respect to the $L^2$-metric $\langle\cdot,\cdot\rangle_{\hyp}$, which is given by equation \eqref{L^2-metric}. Then, for $z,w\in \xbar$, the Bergman kernel associated to the line bundle $\cL^{k}$, is given by the following formula
\begin{align*}
\bkl(z,w):=\sum_{j=1}^{n_k}\omega_j(z)\overline{\omega_j(w)}.
\end{align*}

The Bergman kernel $\bkl$ is the distribution kernel of the orthogonal projection onto $H^0\big(\xbar,\cL^{k}\big)$. So by Reisz representation theorem, the definition of the Bergman kernel is independent of the choice of any orthonormal basis for $H^0\big(\xbar,\cL^{k}\big)$. Furthermore, the Bergman kernel $\bkl$ is a holomorphic differential form in the $z$-variable, for $z\in\xbar\backslash\lbrace i\infty\rbrace$, and an anti-holomorphic differential form in the $w$-variable,  for $w\in\xbar\backslash\lbrace i\infty\rbrace$.

\vspace{0.1cm}
The point-wise hyperbolic metric on $H^0\big(\xbar, \cL^k\big)$ induces a point-wise hyperbolic metric on $\bkl$ along the diagonal, which at $z\in \xbar$,  
is given by the following formula
\begin{align*}
\|\bkl(z)\|_{\hyp}=\sum_{j=1}^{n_k}\| \omega_j(z) \|_{\hyp}^2.
\end{align*}   

We now define the Bergman metric associated to $\cL^k$. For $k\geq 1$, and $z\in\xbar \backslash\lbrace i\infty\rbrace$, the Bergman metric associated to the line bundle $\cL^k$, is given by the following formula
\begin{align}\label{muberk}
\muberkbarx(z):=-\frac{i}{2\pi}\partial_{z}\partial_{\overline{z}}\log\|\bkl(z)\|_{\hyp}
\end{align}

We now recall the Bergman kernel associated to $\Sk$. For $z=x+iy\in X$ (identifying $X$ with its universal cover $\H$), and $f\in \Sk$, the point-wise 
Petersson metric $\|\cdot\|_{\mathrm{pet}}$ on $\Sk$, is given by the following formula
\begin{align*}
\| f(z) \|_{\mathrm{pet}}^{2}:= y^{2k}\big|f(z)\big|^{2},
\end{align*}
which in turn induces an $L^2$-metric on $\Sk$, namely the Petersson inner-product on $\Sk$. For $f,g\in\Sk$, the Petersson inner-product $\langle f,g\rangle_{\mathrm{pet}}$ on $\Sk$, is given by the following formula
\begin{align}\label{petip}
\big\langle f,g \big\rangle_{\mathrm{pet}}:=\int_{\cF_{\Gamma}}y^{2k} f(z)\overline{g(z)}\hypbarx(z).
\end{align} 

\vspace{0.1cm}
Let $\big\lbrace{f_1,\cdots, f_{n_k}\big\rbrace}$ be an orthonormal basis of $\Sk$, with respect to the Petersson inner-product. Then for $z,w\in \mathbb{H}$, the Bergman kernel $\bk(z,w)$ associated to the complex vector space $\Sk$, is given by the following formula
\begin{align*}
\bk(z,w):=\sum_{j=1}^{n_k}f_j(z)\overline{f_j(w)}.
\end{align*}

The Bergman kernel $\bk$ is a distribution kernel for the Hilbert space $\Sk$, and hence, by Riesz representation theorem, the definition of the Bergman kernel is independent of the choice of orthonormal basis for $\Sk$. Furthermore, the Bergman kernel $\bk$ is a holomorphic cusp form in the $z$-variable, and an anti-holomorphic cusp form in the $w$-variable. When $z=w$, for brevity of notation, we denote $\bk(z,w)$ by $\bk(z)$. 

\vspace{0.1cm}
The point-wise Petersson metric on $\Sk$ induces a point-wise Petersson metric on $\bk$ along the diagonal, which at $z=x+iy\in X$ (identifying $X$ with its universal cover $\H$), is given by the following formula
\begin{align}\label{bkpet}
\|\bk(z)\|_{\mathrm{pet}}:=\sum_{j=1}^{n_k}\| f_j(z)\|_{\mathrm{pet}}^2=\sum_{j=1}^{n_k}y^{2k}\big| f_j(z)\big|^2=y^{2k}\bk(z).
\end{align} 

We now state an alternative description of the Bergman kernel $\bk$, which for $k\geq 2$, and $z,w\in\H$, is given by the following Poincar\'e series, (see Proposition 1.3 on page 77 in \cite{Fr})
\begin{align}\label{bkseries}
\bk(z,w)=\frac{(2k-1)(2i)^{2k}}{4\pi}\sum_{\gamma\in \G}\frac{1}{(z-\overline{\gamma w})^{2k}}\cdot \frac{1}{\big(\overline{cw+d}\big)^{2k}},\,\,\,\mathrm{where} \,\,\gamma=\big(\begin{smallmatrix} a & b\\ c & d \end{smallmatrix}\big) \in \G.
\end{align}
 
 For $k\geq 1$, and $z=x+iy\in\xbar$ (identifying $X$ with its universal cover $\H$), applying isometry \eqref{isometry}, using the fact that
\begin{align*}
\Im\big(\gamma(z)\big)=\frac{y}{\big|\overline{cz+d}\big|^{2}},\,\,\mathrm{where} \,\,\gamma=\big(\begin{smallmatrix} a & b\\ c & d \end{smallmatrix}\big) \in \G,
\end{align*}
and combining equations \eqref{isometry}, \eqref{bkpet}, \eqref{bkseries}, and \eqref{dh}, we have the following equality of Bergman kernels and their estimate
\begin{align}\label{bkequality}
\|\bkl(z)\|_{\hyp}=\|\bk(z)\|_{\mathrm{pet}}=\frac{2k-1}{4\pi}\sum_{\gamma\in \G}\frac{(4y)^{k}}{\big(z-\overline{\gamma z}\big)^{2k}}\cdot\frac{y^{k}}{\big(\overline{cz+d}\big)^{2k}}\leq\notag\\[0.12cm]
\frac{2k-1}{4\pi}\sum_{\gamma\in \G}\frac{\big|4y\,\Im(\gamma z)\big|^{k}}{\big|z-\gamma \overline{z}\big|^{2k}}=\frac{2k-1}{4\pi} \sum_{\gamma\in \G}\frac{1}{\cosh^{2k}\big(\dh(z,\gamma z)\slash 2\big)}.
\end{align} 
\subsection{Auxiliary estimates}\label{subsec-2.2}
In this section, we compute $n_k$, the dimension of $\Sk$, and for $z\in \cF_{\G,k}$, employing similar techniques as in \cite{am} and \cite{am2}, we compute estimates of the Bergman kernel $\bk(z)$. 

\vspace{0.1cm}
In the following proposition, we compute $n_k=\mathrm{dim}_{\C}H^{0}\big(\xbar,\cL^{k}\big)$, which is equal to the dimension of $\Sk$. The following result is well known, and we refer to \cite{miyake}, for further details. However, for the sake of completion, we present the main ideas behind the proof. 

\vspace{0.1cm}
\begin{prop}\label{prop1}
With notation as above, for $k\geq2$, we have 
\begin{align*}
n_{k}:=\mathrm{dim}_{\C}H^{0}\big(\xbar,\cL^{k}\big)=(2k-1)(g-1)+k-1.
\end{align*}
\begin{proof}
For $k\geq2$, since
\begin{align*}
\deg\big(\Ocpt^{\otimes (1-k)}\otimes\cO_{\xbar}\big((1-k)i\infty\big)\big)=(1-k)(2g-2)+(1-k)<0,
\end{align*}
from Serre duality, we have 
\begin{align*}
\dim_{\C} H^1\big(\xbar, \cL^k\big)=\dim_{\C} H^0\big(\xbar, \Ocpt\otimes (\cL^k)^{\ast}\big)
=\dim_{\C} H^0\big(\xbar, \Ocpt^{\otimes (1-k)}\otimes\cO_{\xbar}\big(1-k)i\infty\big)=0.
\end{align*}

Therefore, applying the Riemann-Roch theorem, we derive
\begin{align*}
n_k=\dim_{\C} H^0\big(\xbar, \cL^k\big)=\deg \big(\cL^k\big)+(1-g)
=\deg\big( \Ocptk\otimes \cO_{\xbar}\big((k-1)i\infty\big)\big)+(1-g)=\notag\\[0.11cm]k(2g-2)+(k-1)+(1-g)=(2k-1)(g-1)+k-1,
\end{align*}  
which completes the proof of the proposition.
\end{proof}
\end{prop}

\vspace{0.1cm}
In the following proposition, combining arguments from Main Theorem in \cite{am2}, and Theorem~6.1 in \cite{fjk2}, we derive an estimate for the Bergman kernel. 
For the benefit of the reader, we reproduce the arguments here. 

\vspace{0.1cm}
\begin{prop}\label{prop3}
With notation as above, for $k\geq 3$, and $z\in\xbar$, we have the following estimate
\begin{align}\label{est1:prop3}
\|\bkl(z)\|_{\hyp}\ll_{\G} k+ \sqrt{k}y;
\end{align}
furthermore, we have 
\begin{align}\label{est2:prop3}
\sup_{z\in\xbar}\|\bkl(z)\|_{\hyp}=O_{\G}\big(k^{3/2}\big),
\end{align}
where the implied constants in estimates \eqref{est1:prop3} and \eqref{est2:prop3}, depend only on the Fuchsian subgroup $\G$.
\begin{proof}
For $k\geq 3$, and $z\in\xbar\backslash\lbrace i \infty\rbrace$, from equation \eqref{bkequality}, we have 
\begin{align}\label{prop3-eqn1} 
&\|\bkl(z)\|_{\hyp}=\|\bk(z)\|_{\mathrm{pet}}=\frac{2k-1}{4\pi}+\notag\\[0.12cm]&\frac{2k-1}{4\pi}\sum_{\gamma\in \G_{i\infty}\backslash\lbrace\mathrm{Id}\rbrace}\frac{(4y)^{k}}{\big(z-\overline{\gamma z}\big)^{2k}}\cdot\frac{y^{k}}{\big(\overline{cz+d}\big)^{2k}}+\frac{2k-1}{4\pi}\sum_{\gamma\in \G\backslash \G_{i\infty}}\frac{(4y)^{k}}{\big(z-\overline{\gamma z}\big)^{2k}}\cdot\frac{y^{k}}{\big(\overline{cz+d}\big)^{2k}},
\end{align}
where $\mathrm{Id}$ denotes the identity matrix. 

\vspace{0.1cm}
Substituting $\delta=0$ in equations (1.6)  and (1.7 ) in \cite{am2}, and slightly refining the arguments, we arrive at the following estimates
\begin{align}\label{prop3-eqn2}
\frac{2k-1}{4\pi}\Bigg|\sum_{\gamma\in \G_{i\infty}\backslash\lbrace\mathrm{Id}\rbrace}\frac{(4y)^{k}}{\big(z-\overline{\gamma z}\big)^{2k}}\cdot\frac{y^{k}}{\big(\overline{cz+d}\big)^{2k}}\Bigg|
\leq  \frac{y(2k-1)}{\sqrt{\pi}}\cdot \frac{\Gamma(k-1\slash 2)}{\Gamma(k)};\notag\\[0.12cm]\frac{2k-1}{4\pi}\Bigg|
\sum_{\gamma\in \G\backslash \G_{i\infty}}\frac{(4y)^{k}}{\big(z-\overline{\gamma z}\big)^{2k}}\cdot\frac{y^{k}}{\big(\overline{cz+d}\big)^{2k}}\Bigg|
\leq \frac{2k-1}{4\pi}\bigg(\frac{16}{\cosh^{2k-4}(\rx\slash 4)}+\frac{8}{\cosh^{2k-3}(\rx\slash 2 )}\bigg)+\notag\\[0.12cm]
 \frac{2k-1}{2\pi \sinh^2(\rx\slash 4)}\bigg( \frac{1}{\cosh^{2k-3}(\rx\slash 2)}+\frac{1}{\cosh^{2k-4}(\rx\slash 2)}\bigg).
\end{align}
Combining estimates \eqref{prop3-eqn1} and \eqref{prop3-eqn2} with the observation 
\begin{align*}
\frac{\Gamma(k-1\slash 2)}{\Gamma(k)}=O\bigg(\frac{1}{\sqrt{k}}\bigg),
\end{align*}
completes the proof of estimate \eqref{est1:prop3}.

We now apply the arguments from the proof of Theorem~6.1 in \cite{fjk2}. For a given $\varepsilon>0$, set 
\begin{align*}
U_{\G,\varepsilon}:=\big\lbrace z=x+iy\in\cF \big|\,y\geq 1/\varepsilon\big\rbrace\,\,\mathrm{and}\,\,\cF_{\G,\varepsilon}:=\cF\backslash U_{\G,\varepsilon}.
\end{align*}
From estimate \eqref{est1:prop3}, it is clear that
\begin{align}\label{prop3-eqn4}
\sup_{z\in \cF_{\G,\varepsilon}}\|\bk(z)\|_{\mathrm{pet}}=O_{\G}\big(k/\sqrt{\varepsilon}\big). 
\end{align}
Writing
\begin{align}\label{prop3-eqn5}
&\|\bk(z)\|_{\mathrm{pet}}=g(z)h(z),\notag\\[0.1cm]
&\mathrm{where}\,\,g(z):=y^{2k}|q(z)|^{2},\,\,\mathrm{and}\,\,h(z)=y^{2k}\sum_{j=1}^{n_{k}}\frac{|f_{j}(z)|^{2}}{|q(z)|^{2}},
\end{align}
and observing that the function $h(z)$ is a subharmonic function, we deduce that
\begin{align}\label{prop3-eqn6}
\sup_{z\in U_{\G,\varepsilon}}g(z)=\sup_{z\in \partial U_{\G,\varepsilon}}g(z),
\end{align}
where $ \partial U_{\G,\varepsilon}$ denotes the boundary of the cuspidal neighborhood $U_{\G,\varepsilon}$. Furthermore, from elementary calculus, we 
derive that 
\begin{align}\label{prop3-eqn7}
\sup_{z\in \mathbb{H}}h(z)=y^{2k}|q(z)|^{2}=y^{k}e^{-4\pi y}|_{y=k/2\pi}.
\end{align}
Combining equation \eqref{prop3-eqn5} with observations \eqref{prop3-eqn4}, \eqref{prop3-eqn6}, and \eqref{prop3-eqn7}, we infer that
\begin{align*}
\sup_{z\in\xbar}\|\bk(z)\|_{\mathrm{pet}}=\sup_{z\in\partial U_{\G, k/2\pi}}\|\bk(z)\|_{\mathrm{pet}}=O_{\G}(k^{3/2}).
\end{align*}
which completes the proof of the proposition.
\end{proof}
\end{prop}

\vspace{0.1cm} 
In the following proposition, we ascertain the asymptotics of the Bergman kernel in the neighborhood of a cusp, for large values of $k$. Consequently, we derive a lower bound for the Bergman kernel.  

\vspace{0.1cm} 
\begin{prop}\label{prop4}
With notation as above, for $k\gg1$, and $z\in\xbar$, we have the following estimate
\begin{align}\label{est:prop4}
\|\bkl(z)\|_{\mathrm{hyp}} \gg_{\G}k,
\end{align}
where the implied constant depends only on the Fuchsian subgroup.
\begin{proof}
For $k\gg1$, and a fixed $\varepsilon>0$, and  $z\in\cF_{\G,\varepsilon}$, with notation and arguments as in the proof of Proposition \ref{prop3}, we derive
\begin{align}\label{prop4-eqn1}
\|\bkl(z)\|_{\mathrm{hyp}} =\|\bk(z)\|_{\mathrm{pet}}=\frac{2k-1}{4\pi}-\Bigg|\frac{2k-1}{4\pi}\sum_{\gamma\in \G\backslash\lbrace\mathrm{Id}\rbrace}\frac{(4y)^{k}}{\big(z-\overline{\gamma z}\big)^{2k}}\cdot\frac{y^{k}}{\big(\overline{cz+d}\big)^{2k}}\bigg|\geq\notag\\[0.1cm]
\frac{2k-1}{4\pi}-\frac{2k-1}{4\pi}\bigg(\frac{16}{\cosh^{2k-4}(r_{\G,\varepsilon}\slash 4)}+\frac{8}{\cosh^{2k-3}(r_{\G,\varepsilon}\slash 2 )}\bigg)-\notag\\[0.12cm]
 \frac{2k-1}{2\pi \sinh^2(r_{\G,\varepsilon}\slash 4)}\bigg( \frac{1}{\cosh^{2k-3}(r_{\G,\varepsilon}\slash 2)}+\frac{1}{\cosh^{2k-4}(r_{\G,\varepsilon}\slash 2)}\bigg)
 \gg_{\G}k,\notag\\[0.12cm]\mathrm{where}\,\,r_{\G,\varepsilon}:=\big\lbrace\dh(z,\gamma z)\big|\, z\in\cF_{\G,\varepsilon},\,\mathrm{and}\,\gamma\in\G\backslash\lbrace\mathrm{Id}\rbrace\big\rbrace.
\end{align}
Furthermore, for a fixed $z=x+iy\in U_{\G,\varepsilon}$, as $k$ approaches $\infty$, observe that absolute value of each of the term in the series
\begin{align*}
\frac{2k-1}{4\pi}\sum_{\gamma\in \G_{i\infty}\backslash\lbrace\mathrm{Id}\rbrace}\frac{(2y)^{2k}}{\big(z-\overline{\gamma z}\big)^{2k}}=
\frac{2k-1}{4\pi}\sum_{n\in\Z\backslash\lbrace 0\rbrace}\frac{(2y)^{2k}}{\big(n+2iy\big)^{2k}}
\end{align*}
approach zero.

Hence, for $k\gg1$, and $z\in U_{\G,\varepsilon}$, we deduce that
\begin{align}\label{prop4-eqn2}
&\|\bk(z)\|_{\mathrm{pet}}=\notag\\[0.1cm]&\frac{2k-1}{4\pi}-\frac{2k-1}{4\pi}\Bigg|\sum_{\gamma\in\G\backslash\G_{i\infty}}\frac{(4y)^{k}}{\big(z-\overline{\gamma z}\big)^{2k}}\cdot\frac{y^{k}}{\big(\overline{cz+d}\big)^{2k}}\Bigg|-\frac{2k-1}{4\pi}\Bigg|\sum_{n\in\Z\backslash\lbrace 0\rbrace}
\frac{(2y)^{2k}}{\big(n+2iy\big)^{2k}}\Bigg|=\notag\\[0.1cm]&\frac{2k-1}{4\pi}-\frac{2k-1}{4\pi}\bigg(\frac{16}{\cosh^{2k-4}(r_{\G}\slash 4)}+\frac{8}{\cosh^{2k-3}(r_{\G}\slash 2 )}\bigg)-\notag\\[0.12cm]
&\frac{2k-1}{2\pi \sinh^2(r_{\G}\slash 4)}\bigg( \frac{1}{\cosh^{2k-3}(r_{\G}\slash 2)}+\frac{1}{\cosh^{2k-4}(r_{\G}\slash 2)}\bigg)
-\frac{2k-1}{4\pi}\Bigg|\sum_{n\in\Z\backslash\lbrace 0\rbrace}\frac{(2y)^{2k}}{\big(n+2iy\big)^{2k}}\Bigg|\gg_{\G}k.
\end{align}

Combining estimates \eqref{prop4-eqn1} and \eqref{prop4-eqn2}, completes the proof of the proposition.
\end{proof}
\end{prop}
\subsection{Estimates of the Bergman metric}\label{subsec-2.3}
In this subsection, extending techniques from \cite{ab}, we prove Theorem \ref{mainthm1}.  

\vspace{0.1cm}
For $k\geq 1$, and $z=x+iy\in\xbar\backslash\lbrace i\infty\rbrace$ (identifying $X$ with its universal cover $\H$), from the definition of the Bergman metric $\muberkbarx(z)$, which is as defined in equation \eqref{muberk}, and combining it with equations \eqref{bkpet}, and \eqref{bkequality}, we arrive at the following equality 
\begin{align}\label{mukber1}
\muberkbarx(z)=-\frac{i}{2\pi}\partial_{z}\partial_{\overline{z}}\log \|\bkxbar(z)\|_{\hyp}=-\frac{i}{2\pi}\partial_{z}\partial_{\overline{z}}\log (y)^{2k}-
\frac{i}{2\pi}\partial_{z}\partial_{\overline{z}}\log \bk(z).
\end{align}

In the following proposition, we relate the Bergman metric and the hyperbolic metric.

\vspace{0.1cm}
\begin{lem}\label{lem4}
With notation as above, for $k\geq 1$, and $z\in\xbar\backslash\lbrace i\infty\rbrace$, we have
\begin{align*}
\muberkbarx(z)=\frac{k}{2\pi}\hypx(z)+\frac{y^{2}}{\pi}\bigg(\displaystyle\frac{\displaystyle\frac{\partial \bk(z)}{\partial z}\frac{\partial \bk(z)}{\partial \overline{z}}}{\big(\bk(z)\big)^{2}}-\displaystyle\frac{\displaystyle\frac{\partial^{2}\bk(z)}{\partial z\partial \overline{z}}}{\bk(z)}\bigg)\hypx(z).
\end{align*}
\begin{proof}
The proof of the proposition follows from combining equation \eqref{mukber1}, and arguments from Proposition 2.1 from \cite{ab}.
\end{proof}
\end{lem}

\vspace{0.1cm}
\begin{lem}\label{lem7}
With notation as above, for $k\geq 3$, and $z=x+ iy\in\H$, we have the following estimate 
\begin{align}\label{est1:lem7}
\bigg|\frac{\partial^{2} \bk(z)}{\partial z\partial \overline{z}} \bigg|\ll_{\G}\frac{k^{3}}{y^{2k+2}}+ \frac{k^{5/2}}{y^{2k+1}};
\end{align}
furthermore, we have
\begin{align}\label{est2:lem7}
\sup_{z\in\xbar}y^{2k+2}\bigg|\frac{\partial^{2} \bk(z)}{\partial z\partial \overline{z}} \bigg|=O_{\G}\big(k^{7/2}\big),
\end{align}
where the implied constants in estimates \eqref{est1:lem7} and \eqref{est2:lem7}, depend only on the Fuchsian subgroup $\G$.
\begin{proof}
Estimate \eqref{est1:lem7} follows from combining from Lemma 2.4 from \cite{ab}, and arguments from the proof of Proposition \ref{prop3}. 

Observing that 
\begin{align*}
\bigg|\frac{\partial^{2} \bk(z)}{\partial z\partial \overline{z}} \bigg|=\bigg|\sum_{j=1}^{n_{k}}\frac{\partial f_{j}(z)}{\partial z}\frac{\partial \overline{f_{j}(z)}}{\partial \overline{z}}\bigg|=\sum_{j=1}^{n_{k}}\bigg|\frac{\partial f_{j}(z)}{\partial z}\bigg|^{2}
\end{align*}
is a sum of subharmonic functions, and hence, a subharmonic function. 

So employing the same arguments as in  the proof of Proposition \ref{prop4}, we can conclude that 
\begin{align*}
\sup_{z\in\xbar}y^{2k+2}\bigg|\frac{\partial^{2} \bk(z)}{\partial z\partial \overline{z}} \bigg|=\sup_{\substack{z=x+iy\in\xbar\\y=(k+1)/\pi}}y^{2k+2}\bigg|\frac{\partial^{2} \bk(z)}{\partial z\partial \overline{z}} \bigg|.
\end{align*}
Combining estimate \eqref{est2:lem7} with the above observation, completes the proof of the lemma.
\end{proof}
\end{lem}

\vspace{0.1cm}
\begin{lem}\label{lem5}
With notation as above, we have the following estimate
\begin{align}\label{est:lem5}
\sup_{z\in\xbar}y^{4k+2}\bigg|\frac{\partial \bk(z)}{\partial z}\frac{\partial \bk(z)}{\partial \overline{z}} \bigg|=O_{\G}(k^{5}),
\end{align}
where the implied constant depend only on the Fuchsian subgroup $\G$.
\begin{proof}
Applying Cauchy-Schwartz inequality, we observe that
\begin{align*}
&\bigg|\frac{\partial \bk(z)}{\partial z}\frac{\partial \bk(z)}{\partial \overline{z}} \bigg|=\bigg|\frac{\partial \bk(z)}{\partial z} \bigg|^{2}= \bigg|\sum_{j=1}^{n_{k}}\frac{\partial f_{j}(z)}{\partial z}\overline{f_{j}(z)}\bigg|^{2}\leq\notag\\[0.12cm] &\bigg(
\sum_{j=1}^{n_{k}}\big|f_{j}(z)\big|^{2}\bigg)\bigg(\sum_{j^{\prime}=1}^{n_{k}}\bigg|\frac{\partial f_{j^{\prime}}(z)}{\partial z}\bigg|^{2}\bigg)=\bk(z) \bigg|\frac{\partial^{2} \bk(z)}{\partial z\partial \overline{z}} \bigg|.
\end{align*}
Combining estimates \eqref{est2:prop3} and \eqref{est2:lem7} with the above observation, completes the proof of the lemma. 
\end{proof}
\end{lem}

\vspace{0.1cm}
Applying the above results, we now estimate the Bergman metric.

\vspace{0.1cm}
\begin{thm}\label{thm10}
With notation as above, for $k\gg 1$, we have the following estimate
\begin{align*}
\sup_{z\in \xbar}\bigg|\frac{\muberkbarx(z)}{\hypbarx( z)} \bigg|=O_{\G}(k^{3}),
\end{align*}
where the implied constant depends only on the Fuchsian subgroup $\G$.
\begin{proof}
For $k\geq3$, and $z\in\xbar\backslash\lbrace i\infty\rbrace$, from Lemma \ref{lem4}, we have
\begin{align}\label{prop8-eqn1}
\bigg|\frac{\muberkbarx(z)}{\hypbarx( z)} \bigg|\leq\frac{k}{2\pi}+\frac{y^{2}}{\pi}\Bigg|\displaystyle\frac{\displaystyle\frac{\partial \bk(z)}{\partial z}\frac{\partial \bk(z)}{\partial \overline{z}}}{\big(\bk(z)\big)^{2}}\Bigg|+\frac{y^{2}}{\pi}\Bigg|\displaystyle\frac{\displaystyle\frac{\partial^{2}\bk(z)}{\partial z\partial \overline{z}}}{\bk(z)}\Bigg|.
\end{align}
We now estimate the second term on the right-hand side of the above inequality. For $z\in \xbar$, from estimate \eqref{est:lem5}, we derive
\begin{align}\label{prop8-eqn2}
\frac{y^{2}}{\pi}\Bigg|\frac{\displaystyle\frac{\partial \bk(z)}{\partial z}\frac{\partial \bk(z)}{\partial \overline{z}}}{\big(\bk(z)\big)^{2}}\Bigg|
=\frac{y^{4k+2}\bigg|\displaystyle\frac{\partial \bk(z)}{\partial z}\frac{\partial \bk(z)}{\partial \overline{z}}\bigg|}{\pi\|\bk(z)\|_{\mathrm{pet}}^{2}}\ll_{\G}\frac{k^{5}}{\|\bk(z)\|_{\mathrm{pet}}^{2}}.
\end{align}
Furthermore, for $k\gg 1$, from estimate \eqref{est:prop4}, we have
\begin{align}\label{prop8-eqn4}
\|\bk(z)\|_{\mathrm{pet}}=\|\bkl(z)\|_{\mathrm{hyp}}\gg_{\G} k,
\end{align}
Combining estimates \eqref{prop8-eqn2} and \eqref{prop8-eqn4}, we arrive at the following estimate, for the second term on the right-hand side of inequality \eqref{prop8-eqn1}
\begin{align}\label{prop8-eqn5}
\frac{y^{2}}{\pi}\Bigg|\frac{\displaystyle\frac{\partial \bk(z)}{\partial z}\frac{\partial \bk(z)}{\partial \overline{z}}}{\big(\bk(z)\big)^{2}}\Bigg|=O_{\G}(k^{3}).
\end{align}
We now estimate the third term on the right-hand side of inequality \eqref{prop8-eqn1}. For $z\in \xbar$, combining estimates \eqref{est2:lem7} and \eqref{prop8-eqn4}, we arrive at the following estimate, for the third term on the right-hand side of inequality \eqref{prop8-eqn1}
\begin{align}\label{prop8-eqn3}
\frac{y^{2}}{\pi}\Bigg|\displaystyle\frac{\displaystyle\frac{\partial^{2}\bk(z)}{\partial z\partial \overline{z}}}{\bk(z)}\Bigg|=\frac{y^{2k+2}\bigg|\displaystyle\frac{\partial^{2}\bk(z)}{\partial z\partial \overline{z}}\bigg|}{\pi\|\bk(z)\|_{\mathrm{pet}}}=O_{\G}(k^{5/2}).
\end{align}
Combining equation \eqref{prop8-eqn1} with estimates \eqref{prop8-eqn5} and \eqref{prop8-eqn3}, completes the proof of the theorem. 
\end{proof}
\end{thm}
\section{Comparison of K\"ahler metrics on symmetric product of Riemann surfaces}
\subsection{Symmetric product of noncompact Riemann surfaces}\label{subsec-3.1}
We now set up the notation, and recall the back ground material to prove Main Theorem \ref{mainthm2}. However, the notation and the requisite back ground details, remain the same, as in the compact setting. So, for brevity of the article, we only state the requisite details, and we refer the reader to section 3 in \cite{ab}, for a more detailed description. 

\vspace{0.1cm}
With notation as in section \ref{subsec-2.1}, let $X$ denote a noncompact finite volume hyperbolic Riemann surface of genus $g\geq 2$. Let $\xbar:=X\cup\lbrace i\infty\rbrace$ denote the Satake compactification of $X$, and let $\overline{X}^d:=\xbar\times\cdots\times \overline{X}$ denote the $d$-fold Cartesian product of $\xbar$. For $1\leq i\leq d$, let 
\begin{align*}
p_{i}:\overline{X}^{d}\longrightarrow \overline{X}
\end{align*}
denote the projection to the $i$-th factor. Put
\begin{align*}
\mu_{\overline{X}^d}^{\hyp}:=\sum_{i=1}^{d}p_{i}^{\ast}\mu_{\overline{X}}^{\hyp},
\end{align*}
which defines a K\"ahler metric on $\overline{X}^d$. 

\vspace{0.1cm}
Let $S_d$ denote the group of permutations of the set $\lbrace 1,\ldots,d \rbrace$ of $d$-elements. The permutation group $S_d$ acts on the Cartesian  product $\xbar^d$, and we denote the quotient space $S_d\backslash \xbar^d$ by $\Sym^d(\xbar)$. The symmetric product $\Sym^d(\xbar)$ is 
an irreducible smooth complex projective variety of complex dimension $d$. 

\vspace{0.1cm}
The K\"ahler metric $\mu_{\overline{X}^d}^{\hyp}$ descends to define a K\"ahler metric on $\Sym^d(\xbar)$, which admits singularities, which we again denote by $\mu_{\overline{X}^d}^{\hyp}$, for brevity of notation. Furthermore, let $\mu_{\overline{X}^d,\mathrm{vol}}^{\hyp}$ denote the volume form associated to the (1,1)-form $\mu_{\overline{X}^d}^{\hyp}$ on $\Sym^d(\xbar)$. 

\vspace{0.2cm}
\paragraph{\bf{Embedding of the symmetric product into the Grassmannian}}
From section \ref{subsec-2.1}, for $k\geq 1$, recall the line bundle
\begin{align*}
\cL^k=\Omega_{\overline{X}}^{\otimes k}\otimes \mathcal{O}_{\overline{X}}\big((k-1)i\infty\big)
\end{align*}
defined on $\overline{X}$. 

\vspace{0.1cm}
Recall from Proposition \ref{prop1}, we know that
\begin{align*}
n_{k}=\mathrm{dim}_{\C}H^{0}\big(\xbar,\cL^{k}\big)=(2k-1)(g-1)+k-1.
\end{align*}

\vspace{0.1cm}
Let $d\geq 1$ be a given integer. Then for for $k\geq 1$, and a divisor $D$ on $\overline{X}$ of degree $d$, let
\begin{align*}
\cL^{k}(-D):=\cL^k\otimes\cO_{\xbar}(-D)=\Omega_{\xbar}^{\otimes k}\otimes \cO_{\xbar}\big((k-1)i\infty-D\big).
\end{align*}

Using Riemann-Roch theorem, we now compute 
\begin{align*}
r_{k}:=\dim_{\C}H^0\big(\xbar,\cL^{k}(-D)\big).
\end{align*}

From Serre duality, we have
\begin{align}\label{sd1}
\dim_{\C} H^1(\overline{X},\mathcal{L}^k(-D))=\dim_{\C} H^0\big(\overline{X},\Omega_{\xbar}\otimes\big(\cL^{k}(-D)\big)^{\ast}\big)=\notag\\[0.1cm]
\dim_{\C} H^0\big(\overline{X},\Omega_{\overline{X}}^{\otimes (1-k)}\otimes \mathcal{O}_{\overline{X}}\big(-(k-1)i\infty+D\big)\big).
\end{align}

Furthermore, observe that
\begin{align}\label{sd2}
\deg\big (\Omega_{\overline{X}}^{\otimes (1-k)}\otimes \mathcal{O}_{\overline{X}}\big(-(k-1)i\infty+D\big)\big)=\notag\\[0.1cm]-(k-1)(2g-2)-(k-1)+d=-(k-1)(2g-1)+d.
\end{align}

So, for $(k-1)(2g-1)>d$, combining equations \eqref{sd1} and \eqref{sd2}, we infer that
\begin{align*}
\dim_{\C} H^1(\overline{X},\mathcal{L}^k(-D))=0
\end{align*}

Hence, for $(k-1)(2g-1)>d$, applying Riemann-Roch theorem, we deduce that
\begin{align*}\label{rr1}
r_{k}=\dim_{\C} H^0\big(\overline{X},\mathcal{L}^k(-D)\big)= \deg\big(\mathcal{L}^k(-D)\big)+(1-g)=\notag\\[0.1cm]
k(2g-2)+(k-1)-d+(1-g)=(2k-1)(g-1)+k-1-d=n_k-d.
\end{align*} 

\vspace{0.1cm}
Consider the following short exact sequence of sheaves over $\overline{X}$
\begin{equation}\label{ses1}
0\rightarrow \mathcal{L}^k(-D) \rightarrow \mathcal{L}^k \rightarrow \mathcal{L}^k\mid _{D} \rightarrow 0,
\end{equation} 
which induces the injective homomorphism 
\begin{align}\label{ses2}
H^0\big(\overline{X},\mathcal{L}^k(-D)\big)\hookrightarrow H^0\big(\overline{X},\mathcal{L}^k\big),
\end{align}
which occurs in the long exact sequence of cohomologies associated with equation \eqref{ses1}.  

\vspace{0.1cm}
From equation \eqref{ses2}, for a given $d\geq 1$, and $(k-1)(2g-1)>d$, we have the following map
\begin{align}\label{emb}
\varphi^{k,d}_{\mathcal{L}} : \Sym^d(\overline{X}) \rightarrow &\Gr(r_k,n_k)\notag\\[0.1cm]
\big(z_1,\ldots,z_d\big)\mapsto & H^0(\overline{X},\mathcal{L}^k(-D)\big)\subset H^0\big(\xbar,\mathcal{L}^k\big),  
\end{align}
where $D$ denotes the divisor $x_1+\cdots+x_d$, and $z_j$ denotes the complex coordinate of the point $x_j\in\xbar$, for each $1\leq j \leq d$. 
Furthermore, the map $\varphi^{k,d}_{\mathcal{L}}$ is a holomorphic embedding, (see section 4 \& 5 in \cite{biswas}). 

\vspace{0.2cm}
\paragraph{\bf{Bergman kernel and the Fubini-Study metric}} 
With notation as above, let $D$ be an effective divisor of $\xbar$ of degree $d\geq 1$. For $k\geq1$, let $\bkld$ denote the Bergman kernel associated to the line bundle $\cL^k(-D)$, and $\|\cdot\|_{\hyp}$, the point-wise hyperbolic metric on $H^0\big(X,\cL^k\big)$, induces a metric on  $H^0\big(X,\cL^k(-D)\big)$, which we again denote by $\|\cdot\|_{\hyp}$, for brevity of notation. 

\vspace{0.1cm}
For $k\gg 1$, and a fixed $z\in \xbar\backslash\lbrace i\infty\rbrace$, using Theorem 3.3 from \cite{ma1}, and triangular inequality, we have the following asymptotic relation
\begin{align}\label{ma}
\frac{1}{k}\|\bkld(z)\|_{\hyp}=\frac{1}{k}\|\bkl(z)\|_{\hyp}+O_{\G}\bigg(\frac{1}{k}\bigg),
\end{align}
where the implied constant depends only on the Fuchsian subgroup $\G$. 

\vspace{0.1cm}
With notation as above, let $\mu^{\mathrm{FS}}_{\mathrm{Gr}(r_k,n_k)}$ denote the Fubini-Study metric on $\mathrm{Gr}(r_k,n_k)$. Let 
\begin{align*}
\mu^{\mathrm{FS},k}_{\Sym^d(\xbar)}:=\big(\varphi^{k,d}_{\cL}\big)^{\ast}\mu^{\mathrm{FS}}_{\mathrm{Gr}(r_k,n_k)}
\end{align*}
denote the pull back of the Fubini-Study metric on $\Sym^d(\xbar)$, via the holomorphic embedding $\varphi^{k,d}_{\mathcal{L}}$, as in equation \eqref{emb}. Furthermore, let 
$\mu^{\mathrm{FS},k}_{\Sym^d(\xbar),\mathrm{vol}}$ denote the volume form associated to the (1,1)-form $\mu^{\mathrm{FS},k}_{\Sym^d(\xbar)}$.

\vspace{0.1cm}
With notation as above, for $z:=(z_1,\ldots,z_d)\in\Sym^d(\xbar)$, from Proposition 3.1 from \cite{ab}, we have the following relation
\begin{align}\label{ab1}
\mu^{\mathrm{FS},k}_{\Sym^d(\xbar)}(z)=-\frac{i}{2\pi}\sum_{j=1}^{d}\partial_{z_j}\partial_{\overline{z}_{j}}\log\|\bkld(z_{j})\|_{\mathrm{hyp}}.
\end{align}
\subsection{Proof of Main Theorem \ref{mainthm2}}\label{subsec-3.2}
We now prove Main Theorem \ref{mainthm2} in this section. 

\vspace{0.1cm}
\begin{thm}\label{thm11}
With notation as above, for a given $d\geq 1$, and $k\gg1$, and for $z\in\Sym^d(\xbar)$, we have the following estimate
\begin{align*}
\bigg|\frac{\mu^{\mathrm{FS},k}_{\Sym^d(\xbar)}(z)}{\mu_{\overline{X}^d}^{\hyp}(z)} \bigg|=O_{\G}(k^{3d}),
\end{align*}
where the implied constant depends only on the Fuchsian subgroup $\G$.
\begin{proof}
The proof of the theorem follows from combining equation \eqref{ab1} with estimate \eqref{ma}, and Theorem \ref{thm10}.
\end{proof}
\end{thm}
\section*{Acknowledgements}
Both the authors thank Prof. Biswas and Prof. Nagaraj, for many interesting mathematical discussion, which culminated in the completion of the article. 
Both the authors extend their gratitude to Dr. Aprameyan Parthsarathy for his suggestions and  comments, regarding the article. 


\begin{thebibliography}{AAAAAA}
\bibitem{arbarello} E. Arbarello, M. Cornalba, P. A. Griffiths and J. Harris, 
\newblock Geometry of algebraic curves. Vol. I.,
\newblock Grundlehren der Mathematischen Wissenschaften, 267. Springer- Verlag, New York, 1985.
\bibitem{abms} A. Aryasomayajula, I. Biswas, A. S. Morye, and T. Sengupta,
\newblock On the K\"ahler metrics over $\symxd$,
\newblock  J. Geom. Phys. 110 (2016), 187--194.
\bibitem{ab} A. Aryasomayajula and I. Biswas, 
\newblock Bergman kernel on Riemann surfaces and K\"ahler metric on symmetric products,
\newblock Internat. J. Math. 30 (2019).
\bibitem{am} A. Aryasomayajula and P. Majumder,
\newblock Off-diagonal estimates of the {Bergman} kernel on hyperbolic
{Riemann} surfaces of finite volume,
\newblock Proceedings of Amer. Math. Soc. (2018), 4009--4020.
\bibitem{am2} A. Aryasomayajula and P. Majumder,
\newblock Off-diagonal estimates of the {Bergman} kernel on hyperbolic
{Riemann} surfaces of finite volume-{II},
\newblock Ann. Fac. Sci. Toulouse Math. 29 (2020), 795--804.
\bibitem{auvray} H. Auvray, X. Ma, and G. Marinescu, 
\newblock{\it{Bergman kernels on punctured Riemann surfaces}}, 
\newblock{Math. Ann. 379 (2021), 951--1002}.
\bibitem{auvray1} H. Auvray, X. Ma, and G. Marinescu, 
\newblock{\it{Quotient of Bergman kernels on punctured Riemann surfaces}},
\newblock{Math. Z. 301 (2022), 2339--2367}.
\bibitem{biswas} I. Biswas and N. M. Rom\~ao, 
\newblock Moduli of vortices and Grassmann manifolds,
\newblock Comm. Math. Phys. 320 (2013), 1--20.
\bibitem{Fr} E. Freitag,
\newblock \textit{Hilbert Modular Forms},
\newblock Springer-Verlag, Berlin, 1990.
\bibitem{fjk2} J.~Friedman, J.~Jorgenson, and J.~Kramer, 
\newblock{\emph{Uniform sup-norm bounds on average for cusp forms of higher weights}},
\newblock{Arbeitstagung Bonn 2013, 127--154, Progr. Math. 319, Birkh\"auser/Springer, Cham, 2016}.
\bibitem{hsiao} C.-Y. Hsiao, and G. Marinescu, 
\newblock{\it{Asymptotics of spectral function of lower energy forms and Bergman kernel of semi-positive and big line bundles}}, 
\newblock{Comm. Anal. Geom. 22 (2014), 1--108}.
\bibitem{hsiao1} C.-Y. Hsiao, G. Marinescu, and H. Wang, 
\newblock{\it{Szeg\"o kernel asymptotics on complete strictly pseudoconvex CR manifolds}},
\newblock{J. Geom. Anal. 32 (2022)}.
\bibitem{kempf} G. Kempf,
\newblock{\it{Toward the inversion of abelian integrals. I}},
\newblock {Ann. of Math. 110 (1979), 243--273}.
\bibitem{kord} Y. A. Kordyukov, X. Ma, and G. Marinescu, 
\newblock{\it{Generalized Bergman kernels on symplectic manifolds of bounded geometry}}, 
\newblock{Comm. Partial Differential Equations 44 (2019), 1037--1071}.
\bibitem{ma1} X. Ma, and G. Marinescu, 
\newblock{\it{Exponential Estimate for the asymptotics of Bergman kernels}},
\newblock{Math. Ann. 362 (2015), 1327--1347}.
\bibitem{ma}  X. Ma, and G. Marinescu,
\newblock Holomorphic Morse inequalities and Bergman kernels,
\newblock Prog. Math. 254, Birkh\"auser Verlag, Basel, Switzerland, 2007.
\bibitem{manton} N.S.Manton,
\newblock One-vortex moduli space and Ricci flow.
\newblock J. Geom. Phys. 58 (2008), 1772--1783.
\bibitem{miyake} T. Miyake, 
\newblock Modular forms.
\newblock Springer Monographs in Mathematics. Springer-Verlag, Berlin, 2006. 
\end{thebibliography}
\end{document}